\documentclass{article}
\usepackage[tbtags]{amsmath}
\numberwithin{equation}{section}
\usepackage{amstext}
\usepackage{amssymb}
\usepackage{amsthm}
\usepackage[misc]{ifsym}
\usepackage{cases}
\usepackage{mathrsfs}
\usepackage{bm}
\usepackage{amsmath}
\allowdisplaybreaks[4]
\usepackage{geometry}
\newtheorem{theorem}{\hspace{2em}Theorem}[section]
\newtheorem*{Proof}{\hspace{2em}Proof}
\newtheorem{lemma}{\hspace{2em}Lemma}[section]

\title{A Note on Additive Complements}
\author{Fang-Yu Ma}
\date{}

\begin{document}
	\maketitle
	
	\textbf{Abstract:} Two infinite sequnences $A$ and $B$ of non-negative integers are called additive complements, if their sum 
	contains all sufficiently large integers. Let $A(x)$ and $B(x)$ be the counting functions of $A$ and $B$. In this paper, we extend the results of Liu and Fang in 2016 and obtain some results on additive complements. For example, we prove that there exsit additive complements  $A$ and $B$ such that $\displaystyle\limsup_{x\rightarrow+\infty}\dfrac{A(x)B(x)}{x}=2$ 
	and $A(x)B(x)-x=1$ for infinitely positive integers $x$.
	
	
	\section{Introdution}
	
	\quad Two infinite sequnences $A$ and $B$ of non-negative integers are called additive complements, if their sum 
	$A+B=\{a+b:a\in A, b\in B\}$
	contains all sufficiently large integers. Let $A(x)$ and $B(x)$ be the counting functions of $A$ and $B$, namely
	$$A(x)=\displaystyle\sum_{a\in A \atop a\leq x} 1 \quad \text{and} \quad B(x)=\displaystyle\sum_{b\in B \atop b\leq x} 1.$$
	
	In 1964, Danzer \cite{danzer1964luberein} conjectured that for additive complements $A$ and $B$, if 
	$$\displaystyle\limsup_{x\rightarrow+\infty}\dfrac{A(x)B(x)}{x}\leq 1,$$
	then 
	$$A(x)B(x)-x \rightarrow +\infty   \quad \text{as}  \quad  x \rightarrow +\infty.$$
	In 1994, S$\acute{a}$rk$\ddot{o}$zy and Szemer$\acute{e}$di  \cite{sarkozy1994onaproblem} proved the conjecture. Fang and Chen 
	\cite{chenyong2010onadditive, fangjin2014onadditive} proved the following result: for additive complements $A$ and $B$, if 
	$$\displaystyle\limsup_{x\rightarrow+\infty}\dfrac{A(x)B(x)}{x}>2 \quad \text{or} \quad \displaystyle\limsup_{x\rightarrow+\infty}\dfrac{A(x)B(x)}{x}< 3-\sqrt{3},$$
	\renewcommand{\thefootnote}{}\footnotetext
	{2020 Mathematics Subject Classification: Primary 11B13, 11B34.
			
	\quad Keywords. Additive complements, Counting functions, Upper limit.}
	
	

	then
	$$A(x)B(x)-x \rightarrow +\infty \quad \text{as}  \quad  x \rightarrow +\infty.$$
	
	On the other hand, Chen and Fang \cite{chenyong2011onadditive} proved the following result.
	
	\textbf{ Theorem A.}
	For any integer $a$ with $a \geq 2$, there exist additive complements $A$ and $B$ such that 
	$$ \displaystyle\limsup_{x\rightarrow+\infty}\dfrac{A(x)B(x)}{x}=\frac{2a+2}{a+2} $$
	and $A(x)B(x)-x=1$ for infinitely positive integers $x$.

	In 2016, Liu and Fang \cite{fangjin2016anoteon} extended the above result and proved the following result.
	
	\textbf{ Theorem B.}
	For any integer $a, b$ with $2\leq a\leq b$, there exist additive complements $A$ and $B$ such that 
	$$ \displaystyle\limsup_{x\rightarrow+\infty}\dfrac{A(x)B(x)}{x}=\frac{2}{\frac{a-1}{ab-1}+1} $$
	and $A(x)B(x)-x=1$ for infinitely positive integers $x$.
	
    In this paper, we extend Theorem B and prove the following results.
    \begin{theorem}
    	There exist additive complements $A$ and $B$ such that 
    	$$\displaystyle\limsup_{x\rightarrow+\infty}\dfrac{A(x)B(x)}{x}=2$$
    	and $A(x)B(x)-x=1$ for infinitely positive integers $x$.
    \end{theorem}
    \begin{theorem}
    	There exist additive complements $A$ and $B$ such that 
    	$$\displaystyle\limsup_{x\rightarrow+\infty}\dfrac{A(x)B(x)}{x}\in (\frac{16}{9},2)\backslash \mathbb{Q}$$
        and $A(x)B(x)-x=1$ for infinitely positive integers $x$.
    \end{theorem}
    Let 
    \begin{equation}
    	\begin{aligned}
    		\mathcal{L}=
    		& \left\{\displaystyle\limsup_{x\rightarrow+\infty}\dfrac{A(x)B(x)}{x}: A , B \ \text{are additive complements,}\ A(x)B(x)-x=1 \ \text{for infinitely} \ x \right\}, \nonumber
    	\end{aligned}
    \end{equation}
    and $\mathcal{L}^{'}$ be the derivative of $\mathcal{L}$, that is, the set of all the clusterpoints of $\mathcal{L}$. By Theorem A and Theorem B, we have $2\in \mathcal{L}^{'}.$
    \begin{theorem}
    	For any integer $a,b$ with $b>a\geq 2$ and $b\geq a+2$, we have $\dfrac{2}{1+\frac{a}{b(a+1)}}\in \mathcal{L}^{'}.$
    \end{theorem}
\section{Proof of Theorem}

    \quad We will use the following lemma and prove it in section $3$.
    
    \begin{lemma}
    	For any sequence 	$\{b_{j}\}$ of positive integers with $b_{0}=1,b_{j}\geq 2,j\geq 1$, there exist additive complements 
    	$A$ and $B$ such that 
    	$$\displaystyle\limsup_{x\rightarrow+\infty}\dfrac{A(x)B(x)}{x}=
    	\displaystyle\limsup_{k\rightarrow+\infty}\frac{2}{1+D_{k}},$$
    	where 
    	\begin{equation}
    		D_{k}=\displaystyle\sum_{i=0}^{k-1} (-1)^{i} (\prod \limits_{j=0}^{i} b_{k-j})^{-1},
    	\end{equation}
    	and $A(x)B(x)-x=1$ for infinitely positive integers $x$.
    \end{lemma}
    
    For formula $(2.1)$, the expansion is 
    $$D_{k}=\frac{1}{b_{k}}+\frac{(-1)}{b_{k}b_{k-1}}+\frac{(-1)^2}{b_{k}b_{k-1}b_{k-2}}+\cdots+\frac{(-1)^{k-1}}{b_{k}b_{k-1}b_{k-2}\cdots b_{1}}.$$
    
    \textbf{Proof of Theorem $1.1$} \quad For any sequence $\{b_{j}\}$ of positive integers with $b_{0}=1,b_{j}\geq 2,j \geq 1$ and $\displaystyle \liminf_{j\rightarrow+\infty}\frac{1}{b_{j}}=0$, by $(2.1)$ we have 
    $\frac{1}{b_{j}}>D_{j}>0$, then $\displaystyle \liminf_{j\rightarrow+\infty}D_{j}=0$. According to Lemma $2.1$, we have 
    $$\displaystyle\limsup_{x\rightarrow+\infty}\dfrac{A(x)B(x)}{x}=
    \displaystyle\limsup_{j\rightarrow+\infty}\frac{2}{1+D_{j}}=2.$$
    This completes the proof.
    
    \textbf{Proof of Theorem $1.2$}\quad For any sequence $\{d_{i}\}$ of positive integers with $d_{i}\in \{a,b\}$ for all $i\geq 1$, $b>a\geq 2, b\geq a+2 $, we construct a new sequence $\{b_{k}\}$ as 
    $$d_{1},c,d_{2},d_{1},c,d_{3},d_{2},d_{1},c,d_{4},d_{3},d_{2},d_{1},c, \cdots$$
    where the positive integer $c>2b$. Let the items that are equal to $c$ in $\{b_{k}\}$ form a subsequence $\{b_{k_{j}}\}$.
    By Lemma $2.1$, we only need to compute $\displaystyle \liminf_{k\rightarrow+\infty}D_{k}$ for the sequence $\{b_{k}\}$.
    
    If $k\neq k_{j}$, then $b_{k}=a$ or $b$. If $b_{k}=a$, then $D_{k}\geq \frac{1}{b_{k}}-\frac{1}{b_{k}b_{k-1}}\geq \frac{1}{a}-\frac{1}{a^{2}}$. By $b \geq a+2$ and $a\geq 2$, we have $\frac{1}{a}-\frac{1}{a^2}\geq \frac{1}{b}> \frac{1}{c}$. If $b_{k}=b$, then $D_{k}\geq \frac{1}{b_{k}}-\frac{1}{b_{k}b_{k-1}}\geq \frac{1}{b}-\frac{1}{ba}=\frac{1}{b}(1-\frac{1}{a})\geq \frac{1}{2b}>\frac{1}{c}$. $D_{k_{j}}\leq \frac{1}{b_{k_{j}}}=\frac{1}{c}$,
    so $\displaystyle \liminf_{k\rightarrow +\infty}D_{k}=\displaystyle \liminf_{j\rightarrow+\infty}D_{k_{j}}$.
    
    $$D_{k_{j}}=\frac{1}{c}+\frac{(-1)}{cd_{1}}+\frac{(-1)^{2}}{cd_{1}d_{2}}+\frac{(-1)^{3}}{cd_{1}d_{2}d_{3}}+\cdots+
    \frac{(-1)^{j}}{cd_{1}d_{2}d_{3}\cdots d_{j}}+\frac{(-1)^{j+1}}{c^{2}d_{1}d_{2}d_{3}\cdots d_{j}}+\cdots,$$
    so
    \begin{equation}
    	\begin{aligned}
    		\bigg| D_{k_{j}}-\left(\frac{1}{c}+\frac{(-1)}{cd_{1}}+\frac{(-1)^{2}}{cd_{1}d_{2}}+\frac{(-1)^{3}}{cd_{1}d_{2}d_{3}}
    		+\cdots+ \frac{(-1)^{j}}{cd_{1}d_{2}d_{3}\cdots d_{j}}\right) \bigg| \leq \frac{1}{c^{2}d_{1}d_{2}d_{3}\cdots d_{j}}.\nonumber
    	\end{aligned}
    \end{equation}
    Let $j\rightarrow +\infty$, then $\frac{1}{c^{2}d_{1}d_{2}d_{3}\cdots d_{j}}\rightarrow 0$ as $d_{i}\geq 2$. Hence
    \begin{equation}
    	\begin{aligned}
    	\displaystyle \liminf_{j\rightarrow+\infty}D_{k_{j}}=\frac{1}{c}\displaystyle\liminf_{j\rightarrow+\infty}\left(1-\frac{1}{d_{1}}+\frac{1}{d_{1}d_{2}}-
    	\frac{1}{d_{1}d_{2}d_{3}}+\cdots+\frac{(-1)^{j}}{d_{1}d_{2}d_{3}\cdots d_{j}}\right).\nonumber
        \end{aligned}
    \end{equation}
    By $d_{i}\geq 2$ we can easily prove the existence of $\displaystyle \lim_{j\rightarrow+\infty} \left(1-\frac{1}{d_{1}}+\frac{1}{d_{1}d_{2}}-
    \frac{1}{d_{1}d_{2}d_{3}}+\cdots+\frac{(-1)^{j}}{d_{1}d_{2}d_{3}\cdots d_{j}}\right)$. Let
    $$\Delta=\displaystyle \lim_{j\rightarrow+\infty}\left(1-\frac{1}{d_{1}}+\frac{1}{d_{1}d_{2}}-
    \frac{1}{d_{1}d_{2}d_{3}}+\cdots+\frac{(-1)^{j}}{d_{1}d_{2}d_{3}\cdots d_{j}}\right).$$
    Obviously, $1>\Delta \geq 1-\frac{1}{d_{1}}\geq \frac{1}{2}$, $c>2b\geq 2(a+2) \geq 8$, so
    $$0\leq \displaystyle\liminf_{j\rightarrow+\infty}D_{k_{j}}=\frac{\Delta}{c}<\frac{1}{8}, \quad \frac{16}{9}<\displaystyle\limsup_{x\rightarrow+\infty}\dfrac{A(x)B(x)}{x} \leq 2.$$
    
    Next we prove that if $\{d_{i}\}$ and $\{d_{i}^{'}\}$ are different, the corresponding $\Delta$ are not equal. We note that $\Delta$ corresponding to $\{d_{i}\}$ is $\Delta(\{d_{i}\})$ and $\Delta$ corresponding to $\{d_{i}^{'}\}$ is $\Delta(\{d_{i}^{'}\})$.
    
    If the first term of $\{d_{i}\}$ and $\{d_{i}^{'}\}$ are not equal, that is, $d_{1}\neq d_{1}^{'}$. Suppose $d_{1}=a,d_{1}^{'}=b$, then
    $$\Delta(\{d_{i}\})=\displaystyle \lim_{j\rightarrow+\infty}\left(1-\frac{1}{d_{1}}+\frac{1}{d_{1}d_{2}}-
    \frac{1}{d_{1}d_{2}d_{3}}+\cdots+\frac{(-1)^{j}}{d_{1}d_{2}d_{3}\cdots d_{j}}\right)< 1-\frac{1}{a}+\frac{1}{ad_{2}}\leq 1-\frac{1}{a}+\frac{1}{a^{2}},$$
    $$\Delta(\{d_{i}^{'}\})=\displaystyle \lim_{j\rightarrow+\infty}\left(1-\frac{1}{d_{1}^{'}}+\frac{1}{d_{1}^{'}d_{2}^{'}}-
    \frac{1}{d_{1}^{'}d_{2}^{'}d_{3}^{'}}+\cdots+\frac{(-1)^{j}}{d_{1}^{'}d_{2}^{'}d_{3}^{'}\cdots d_{j}^{'}}\right)> 1-\frac{1}{b}\geq 1-\frac{1}{a}+\frac{1}{a^{2}},$$
    thus $\Delta(\{d_{i}\})\neq \Delta(\{d_{i}^{'}\})$. 
    
    If the first $m$ terms of $\{d_{i}\}$ and $\{d_{i}^{'}\}$ are the same, but item $(m+1)$th are different. Suppose 	$$d_{m+1}=a,d_{m+1}^{'}=b,\Delta_{m}=1-\frac{1}{d_{1}}+\frac{1}{d_{1}d_{2}}-
    \frac{1}{d_{1}d_{2}d_{3}}+\cdots+\frac{(-1)^{m}}{d_{1}d_{2}d_{3}\cdots d_{m}}.$$
    If $m$ is an even number, then
    \begin{equation}
    	\begin{aligned}
    		\Delta(\{d_{i}\})
    		&= \Delta_{m}+\frac{(-1)^{m}}{d_{1}d_{2}d_{3}\cdots d_{m}}(\frac{-1}{a})+\frac{(-1)^{m}}{d_{1}d_{2}d_{3}\cdots d_{m}}(\frac{(-1)^{2}}{ad_{m+2}})+\frac{(-1)^{m}}{d_{1}d_{2}d_{3}\cdots d_{m}}(\frac{(-1)^{3}}{ad_{m+2}d_{m+3}})+\cdots\\
    		&<  
    		\Delta_{m}+\frac{1}{d_{1}d_{2}d_{3}\cdots d_{m}}(-\frac{1}{a}+\frac{1}{ad_{m+2}})\leq \Delta_{m}+\frac{1}{d_{1}d_{2}d_{3}\cdots d_{m}}(-\frac{1}{a}+\frac{1}{a^{2}}), \nonumber
    	\end{aligned}
    \end{equation}
    
    \begin{equation}
    	\begin{aligned}
    		\Delta(\{d_{i}^{'}\})
    		&= \Delta_{m}+\frac{(-1)^{m}}{d_{1}d_{2}d_{3}\cdots d_{m}}(\frac{-1}{b})+\frac{(-1)^{m}}{d_{1}d_{2}d_{3}\cdots d_{m}}(\frac{(-1)^{2}}{bd_{m+2}^{'}})+
    		\frac{(-1)^{m}}{d_{1}d_{2}d_{3}\cdots d_{m}}(\frac{(-1)^{3}}{bd_{m+2}^{'}d_{m+3}^{'}})+\cdots \\
    		&>  \Delta_{m}+\frac{1}{d_{1}d_{2}d_{3}\cdots d_{m}}(-\frac{1}{b})\geq \Delta_{m}+\frac{1}{d_{1}d_{2}d_{3}\cdots d_{m}}(-\frac{1}{a}+\frac{1}{a^{2}})	,\nonumber
    	\end{aligned}
    \end{equation}
    thus $\Delta(\{d_{i}\})\neq \Delta(\{d_{i}^{'}\})$.	In a similar way, if $m$ is an odd number, we also have $\Delta(\{d_{i}\})\neq \Delta(\{d_{i}^{'}\})$.
    
    According to the construction of $\{d_{i}\}$, it is known that there are uncountable sequences $\{d_{i}\}$.
    So there are uncountable unequal $\Delta$, i.e. there are uncountable unequal real numbers. Thus there must be irrational numbers in 
    these real numbers, i.e. there are irrational numbers in all $\displaystyle\lim_{j\rightarrow+\infty}D_{k_{j}}$. According to Lemma $2.1$, we have $\displaystyle \limsup_{x\rightarrow+\infty}\frac{A(x)B(x)}{x}=\displaystyle \lim_{j\rightarrow+\infty} \frac{2}{1+D_{k_{j}}}$. This completes this proof.
    
    \begin{lemma}
    	For any integer $a,b,l$ with $b>a\geq 2, b \geq a+2$ and $l$ is a positive odd number, there exist additive complements $A$ and $B$ such that 
    	$$\displaystyle\limsup_{x\rightarrow+\infty}\dfrac{A(x)B(x)}{x}=
    	\dfrac{2}{1+\frac{1-\frac{1}{a^{l+1}}}{b(1+\frac{1}{a})(1-\frac{1}{a^{l}b})}},$$
    	and $A(x)B(x)-x=1$ for infinitely positive integers $x$.
    \end{lemma}
    \begin{Proof}
    Let sequence $\{b_{j}\}$ be
    $$\underbrace{a,a,\cdots,a}_{l},b,\underbrace{a,a,\cdots,a}_{l},b,\underbrace{a,a,\cdots,a}_{l},b,\cdots$$
    By $(2.1)$, we have $D_{(l+1)k}<\frac{1}{b_{(l+1)k}}=\frac{1}{b}$ and $D_{t}>\frac{1}{a}-\frac{1}{ac}\geq \frac{1}{a}-\frac{1}{a^{2}}$, where $t \neq (l+1)k$ and $c\in \{a,b\}$. By $b \geq a+2$ and $a\geq 2$, we have $\frac{1}{a}-\frac{1}{a^2}\geq \frac{1}{b}$.
    According to Lemma $2.1$, we have $\displaystyle\limsup_{x \rightarrow+\infty}\dfrac{A(x)B(x)}{x} =\displaystyle\limsup_{k \rightarrow+\infty}\frac{2}{1+D_{(l+1)k}}$, and 
    \begin{equation}
    	\begin{aligned}
    		\displaystyle\lim_{k \rightarrow+\infty}D_{(l+1)k}
    		&= \left(\frac{1}{b}-\frac{1}{ba}+\frac{1}{ba^{2}}+\cdots+\frac{1}{b(-a)^{l}}\right)
    		\left(1+\frac{1}{a^{l}b}+(\frac{1}{a^{l}b})^{2}+(\frac{1}{a^{l}b})^{3}+\cdots\right)\\
    		&=\frac{1-\frac{1}{a^{l+1}}}{b(1+\frac{1}{a})(1-\frac{1}{a^{l}b})}.\nonumber
    	\end{aligned}
    \end{equation}
    Hence $\displaystyle\limsup_{x\rightarrow+\infty}\dfrac{A(x)B(x)}{x}
    =\dfrac{2}{1+\frac{1-\frac{1}{a^{l+1}}}{b(1+\frac{1}{a})(1-\frac{1}{a^{l}b})}}$. This completes the proof.
    \end{Proof}

    \textbf{Proof of Theorem $1.3$} \quad According to Lemma $2.2$ we have $$\dfrac{2}{1+\frac{1-\frac{1}{a^{l+1}}}{b(1+\frac{1}{a})(1-\frac{1}{a^{l}b})}}\in \mathcal{L}.$$ 
    So we have 
    $$\dfrac{2}{1+\frac{a}{b(a+1)}}\in \mathcal{L}^{'} \quad \text{as} \quad l \rightarrow+\infty.$$
    This completes the proof.
    
\section{Proof of Lemma $2.1$}
     \begin{lemma}
    	For any sequence $\{b_{j}\}$ of positive integers with $b_{0}=1,b_{j}\geq2,j\geq1$, every positive integer can be uniquely represented as $\displaystyle \sum_{j=0}^{n} \varepsilon_{j}\prod \limits_{i=0}^{j} b_{i}$, where $\varepsilon_{j}$  are integers satisfying $0\leq \varepsilon_{j} \leq b_{j+1}-1$.
    \end{lemma}
    
    \begin{Proof}
    	Let $a_{0}=1,a_{j}= \prod \limits_{i=0}^{j} b_{i}, j\geq 1$. For any positive integer $m$, there exists an integer $n$ such that $a_{n}\leq m <a_{n+1}$. By division algorithm we have
    	$$m=\varepsilon_{n}a_{n}+r_{n},0 \leq r_{n}<a_{n},0<\varepsilon_{n} \leq b_{n+1}-1.$$
    	If $r_{n}=0$, then $m=\varepsilon_{n}a_{n}$, thus we complete the proof. If $r_{n}\neq 0$, then there exists an integer $n_{1}<n$ such that $a_{n_{1}}\leq r_{n}<a_{n_{1}+1}$. By division algorithm we have
    	$$r_{n}=\varepsilon_{n_{1}}a_{n_{1}}+r_{n_{1}},0\leq r_{n_{1}}<a_{n_{1}},0<\varepsilon_{n_{1}} \leq b_{n_{1}+1}-1.$$
    	If $r_{n_{1}}=0$, then $m=\varepsilon_{n}a_{n}+\varepsilon_{n_{1}}a_{n_{1}}$, thus we complete the proof. If $r_{n_{1}}\neq 0$, we repeat the above operations and in limited steps, we have
    	$$m=\displaystyle \sum_{j=0}^{n} \varepsilon_{j}a_{j}=\displaystyle \sum_{j=0}^{n} \varepsilon_{j}\prod \limits_{i=0}^{j} b_{i},\varepsilon_{j}=0,1,\cdots,b_{j+1}-1,\varepsilon_{n} \neq 0.$$
    	
    	Next we prove the uniqueness. Assume that $m=\displaystyle \sum_{j=0}^{n} \varepsilon_{j}a_{j}=\displaystyle \sum_{i=0}^{k} \varepsilon_{i}^{'} a_{i},\varepsilon_{n} \neq 0,\varepsilon_{k}^{'} \neq 0$. If $n\neq k$,  suppose $n>k$, then $\varepsilon_{n}a_{n}>\displaystyle \sum_{i=0}^{k} \varepsilon_{i}^{'} a_{i}$. Thus $n=k$ and we have 
    	$$\varepsilon_{0}-\varepsilon_{0}^{'}=\displaystyle \sum_{j=1}^{n}(\varepsilon_{j}^{'}-\varepsilon_{j})a_{j}.$$
    	By $b_{1}\Big| \displaystyle \sum_{j=1}^{n}(\varepsilon_{j}^{'}-\varepsilon_{j})a_{j}$, we have $b_{1}\Big| \varepsilon_{0}-\varepsilon_{0}^{'}$, where $\varepsilon_{0},\varepsilon_{0}^{'} \in \{0,1,\cdots,b_{1}-1\}$. Thus $\varepsilon_{0}=\varepsilon_{0}^{'}$ and we have 
    	$$\varepsilon_{1}-\varepsilon_{1}^{'}=\displaystyle \sum_{j=2}^{n}(\varepsilon_{j}^{'}-\varepsilon_{j})(b_{2}b_{3}\cdots b_{j}).$$
    	By $b_{2}\Big| \displaystyle \sum_{j=2}^{n}(\varepsilon_{j}^{'}-\varepsilon_{j})(b_{2}b_{3}\cdots b_{j})$,
    	we have $b_{2}\Big| \varepsilon_{1}-\varepsilon_{1}^{'}$, where $\varepsilon_{1},\varepsilon_{1}^{'} \in \{0,1,\cdots,b_{2}-1\}$. Thus $\varepsilon_{1}=\varepsilon_{1}^{'}$ and we have 
    	$$\varepsilon_{2}-\varepsilon_{2}^{'}=\displaystyle \sum_{j=3}^{n}(\varepsilon_{j}^{'}-\varepsilon_{j})(b_{3}b_{4}\cdots b_{j}).$$
    	We repeat the above operations and in limited steps, we have $\varepsilon_{j}=\varepsilon_{j}^{'},j=0,1,\cdots,n$. This completes the proof.		
    \end{Proof}		
    For any sequence 	$\{b_{j}\}$ of positive integers with $b_{0}=1,b_{j}\geq 2,j\geq 1$, let $a_{0}=1,a_{j}= \prod \limits_{i=0}^{j} b_{i}, j\geq 1$. By Lemma
    $3.1$, we can construct additive complements
    \begin{equation}
    	\begin{aligned}
    		A&=\Big\{\displaystyle \sum_{j=0}^{n} \varepsilon_{2j}a_{2j}, \varepsilon_{2j}=0,1,\cdots,b_{2j+1}-1 \Big\},\\	
    		B&=	\Big\{\displaystyle \sum_{j=0}^{n} \varepsilon_{2j+1}a_{2j+1}, \varepsilon_{2j+1}=0,1,\cdots,b_{2j+2}-1\Big\},	
    	 \end{aligned}
    \end{equation}
    where the summation is a finite sum. Let
    \begin{align}
    y_{k}&=(b_{1}-1)+(b_{3}-1)a_{2}+(b_{5}-1)a_{4}+\cdots+(b_{2k-1}-1)a_{2k-2}+a_{2k},\\
    z_{k}&=(b_{2}-1)a_{1}+(b_{4}-1)a_{3}+(b_{6}-1)a_{5}+\cdots+(b_{2k}-1)a_{2k-1}+a_{2k+1},
    \end{align}
    obviously, $\{y_{k}\}\subset A, \{z_{k}\} \subset B$.
    
    From here to the end, the additive complements $A$ and $B$ are given in $(3.1)$, $D_{k}$ is given in $(2.1)$, $y_{k}$ and $z_{k}$ are given in $(3.2)$ and $(3.3)$.
    
    \begin{lemma}
    	For the additive complements $A$ and $B$, we have
    	$$\frac{A(y_{k})B(y_{k})}{y_{k}} = \frac{2}{1+D_{2k}^{*}},\quad \frac{A(z_{k})B(z_{k})}{z_{k}} = \frac{2}{1+D_{2k+1}^{*}},$$
    	where $D_{2k}^{*}=D_{2k},D_{2k+1}^{*}=D_{2k+1}-\frac{1}{b_{2k+1}b_{2k}b_{2k-1}\cdots b_{2}b_{1}}.$
    \end{lemma}
    
    \begin{Proof}
    	For $y_{k}$, we have $A(y_{k})=2b_{1}b_{3}b_{5}\cdots b_{2k-1},B(y_{k})=b_{2}b_{4}b_{6}\cdots b_{2k}$. Then
    	$$\frac{A(y_{k})B(y_{k})}{y_{k}}=\dfrac{2a_{2k}}{(b_{1}-1)+(b_{3}-1)a_{2}+(b_{5}-1)a_{4}+\cdots+(b_{2k-1}-1)a_{2k-2}+a_{2k}}=\dfrac{2}{1+D_{2k}^{*}},$$    	    	  	  
    	where $D_{2k}^{*}=\displaystyle \sum_{i=0}^{k-1}(b_{2i+1}-1)\frac{a_{2i}}{a_{2k}}=\displaystyle \sum_{i=0}^{2k-1} (-1)^{i} (\prod \limits_{j=0}^{i} b_{2k-j})^{-1}$. 
    	
    	For $z_{k}$, we have $A(z_{k})=b_{1}b_{3}b_{5}\cdots b_{2k-1}b_{2k+1},B(z_{k})=2b_{2}b_{4}b_{6}\cdots b_{2k}$. Then
    	$$\frac{A(z_{k})B(z_{k})}{z_{k}}=\dfrac{2a_{2k+1}}{(b_{2}-1)a_{1}+(b_{4}-1)a_{3}+(b_{6}-1)a_{5}+\cdots+(b_{2k}-1)a_{2k-1}+a_{2k+1}}=\dfrac{2}{1+D_{2k+1}^{*}},$$   	  	
    	where $D_{2k+1}^{*}=\displaystyle \sum_{i=0}^{k-1}(b_{2i+2}-1)\frac{a_{2i+1}}{a_{2k+1}}=\displaystyle \sum_{i=0}^{2k-1}(-1)^{i}(\prod \limits_{j=0}^{i} b_{2k+1-j})^{-1}$. 
    	This completes the proof.
    \end{Proof}	
    
    \begin{lemma}
    	For any integers $a_{1},a_{2},b_{1},b_{2}$ and $u$, if $a_{1}b_{2}-a_{2}b_{1}\geq 0$ and $a_{2}x+b_{2}>0$ for all  $0\leq x \leq u$, then 
    	for all $0\leq x\leq u$ we have
    	$$\frac{a_{1}x+b_{1}}{a_{2}x+b_{2}}\leq \frac{a_{1}u+b_{1}}{a_{2}u+b_{2}}.$$
    \end{lemma}
    
    \begin{lemma}
    For the additive complements $A$ and $B$, if
    $$x=\varepsilon_{0}+\varepsilon_{2}a_{2}+\varepsilon_{4}a_{4}+\cdots+\varepsilon_{2k-2}a_{2k-2}+\varepsilon_{2k}a_{2k} \in A, 0 \leq \varepsilon_{j}\leq b_{j+1}-1,$$
    then
    $$\frac{A(x)B(x)}{x}\leq \frac{A(y_{k})B(y_{k})}{y_{k}}.$$
    \end{lemma}

    \begin{Proof}
    	By noting that $x\in A$, we only need to consider $x \neq y_{k}$, i.e. the following two conditions.
    	
    	\textbf{Case $1$}  $x=\varepsilon_{0}+\varepsilon_{2}a_{2}+\varepsilon_{4}a_{4}+\cdots+\varepsilon_{2k}a_{2k},2 \leq \varepsilon_{2k}\leq b_{2k+1}-1$	.
    	
    	\textbf{Case $2$}  $x=\varepsilon_{0}+\varepsilon_{2}a_{2}+\cdots+\varepsilon_{2m-2}a_{2m-2}+\displaystyle \sum_{i=m}^{k-1}(b_{2i+1}-1)a_{2i}+a_{2k}$, where $\varepsilon_{2m-2}\leq b_{2m-1}-2,1\leq m \leq k($ if $m=k$, then the sum is zero$)$.
    	
    	For Case $1$, we have $A(x)\leq (\varepsilon_{2k}+1)b_{2k-1}b_{2k-3}\cdots b_{3}b_{1},B(x)=b_{2k}b_{2k-2}\cdots b_{2}.$
    	
    	If $3\leq \varepsilon_{2k} \leq b_{2k+1}-1$, then 
    	$$\frac{A(x)B(x)}{x}\leq\frac{(\varepsilon_{2k}+1)a_{2k}}{\varepsilon_{2k}a_{2k}}=1+\frac{1}{\varepsilon_{2k}}\leq \frac{4}{3} \leq
    	\frac{A(y_{k})B(y_{k})}{y_{k}}.$$
    	
    	If $\varepsilon_{2k}=2$, then $x=\varepsilon_{0}+\varepsilon_{2}a_{2}+\varepsilon_{4}a_{4}+\cdots+\varepsilon_{2k-2}a_{2k-2}+2a_{2k}$, where $\varepsilon_{2k-2}\leq b_{2k-1}-1$. Thus $A(x)\leq 2b_{2k-1}b_{2k-3}\cdots b_{3}b_{1}+(\varepsilon_{2k-2}+1)b_{2k-3}\cdots b_{3}b_{1}$.
    	
    	Hence by Lemma $3.3$ $(x=\varepsilon_{2k-2})$ we have
    	\begin{equation}
    		\begin{aligned}
    			\dfrac{A(x)B(x)}{x}
    			&\leq \dfrac{2a_{2k}+(\varepsilon_{2k-2}+1)b_{2k}a_{2k-2}}{2a_{2k}+\varepsilon_{2k-2}a_{2k-2}}\leq \dfrac{2a_{2k}+(b_{2k-1}-1+1)b_{2k}a_{2k-2}}{2a_{2k}+(b_{2k-1}-1)a_{2k-2}} \\
    			&=\dfrac{3}{2+\frac{1}{b_{2k}}-\frac{1}{b_{2k}b_{2k-1}}}
    			\leq \dfrac{2}{1+D_{2k}^{*}}
    			=\frac{A(y_{k})B(y_{k})}{y_{k}},\nonumber
    		\end{aligned}
    	\end{equation}
        where the last inequality is proved by the fact that $0<\frac{1}{b_{2k}}-\frac{1}{b_{2k}b_{2k-1}}<\frac{1}{2}$ and
    	$0<D_{2k}^{*}-(\frac{1}{b_{2k}}-\frac{1}{b_{2k}b_{2k-1}})<\frac{1}{8}$.
    	
    	For Case $2$, if $m\leq k-1$, we have $B(x)=b_{2k}b_{2k-2}\cdots b_{2}$ and $A(x)\leq b_{2k-1}b_{2k-3}\cdots b_{3}b_{1}+(b_{2k-1}-1)b_{2k-3}\cdots b_{3}b_{1}+(b_{2k-3}-1)b_{2k-5}\cdots b_{3}b_{1}+\cdots+(b_{2m+1}-1)b_{2m-1}\cdots b_{3}b_{1}+(\varepsilon_{2m-2}+1)b_{2m-3}\cdots b_{3}b_{1}$. By Lemma $3.3$ $(x=\varepsilon_{2m-2})$ we have
    	\begin{equation}
    		\begin{aligned}
    			\dfrac{A(x)B(x)}{x}
    			&\leq \dfrac{\mathscr{N}}
    			{a_{2k}+\displaystyle \sum_{i=m}^{k-1}(b_{2i+1}-1)a_{2i}+\varepsilon_{2m-2}a_{2m-2}}
    			\leq  \dfrac{\mathscr{N}^{*}}
    			{a_{2k}+\displaystyle \sum_{i=m}^{k-1}(b_{2i+1}-1)a_{2i}+(b_{2m-1}-2)a_{2m-2}}\\
    			&=\dfrac{2a_{2k}-b_{2k}b_{2k-2}b_{2k-4}\cdots b_{2m}a_{2m-2}}{a_{2k}+\displaystyle \sum_{i=m}^{k-1}(b_{2i+1}-1)a_{2i}+(b_{2m-1}-2)a_{2m-2}}
    			=\dfrac{2-(b_{2k-1}b_{2k-3}\cdots b_{2m+1}b_{2m-1})^{-1}}{1+\displaystyle \sum_{i=m}^{k-1}(b_{2i+1}-1)\frac{a_{2i}}{a_{2k}}+(b_{2m-1}-2)\frac{a_{2m-2}}{a_{2k}} }\\
    			&\leq \dfrac{2}{1+D_{2k}^{*}}
    			=\frac{A(y_{k})B(y_{k})}{y_{k}},\nonumber
    		\end{aligned}
    	\end{equation}
    	where $\mathscr{N}=a_{2k}+b_{2k}(b_{2k-1}-1)a_{2k-2}+b_{2k}b_{2k-2}(b_{2k-3}-1)a_{2k-4}+\cdots+b_{2k}b_{2k-2}\cdots b_{2m+2}(b_{2m+1}-1)a_{2m}\\
    	+b_{2k}b_{2k-2}\cdots b_{2m+2}b_{2m}(\varepsilon_{2m-2}+1)a_{2m-2},$\\
    	$\mathscr{N}^{*}=a_{2k}+b_{2k}(b_{2k-1}-1)a_{2k-2}+b_{2k}b_{2k-2}(b_{2k-3}-1)a_{2k-4}+\cdots+b_{2k}b_{2k-2}\cdots b_{2m+2}(b_{2m+1}-1)a_{2m}\\
    	+b_{2k}b_{2k-2}\cdots b_{2m+2}b_{2m}(b_{2m-1}-2+1)a_{2m-2}.$
    	
    	If $m=k$, then
    	$$x=\varepsilon_{0}+\varepsilon_{2}a_{2}+\varepsilon_{4}a_{4}+\cdots+\varepsilon_{2k-2}a_{2k-2}+a_{2k}, \varepsilon_{2k-2}\leq b_{2k-1}-2.$$
    	Let
    	$$x^{'}=\varepsilon_{0}+\varepsilon_{2}a_{2}+\varepsilon_{4}a_{4}+\cdots+(b_{2k-1}-1)a_{2k-2}+a_{2k},$$
    	namely, $\varepsilon_{2k-2}$ in $x$ becomes $b_{2k-1}-1$ in $x^{'}$ and the rest is the same. We have
    	\begin{align}
    	A(x)&= b_{2k-1}b_{2k-3}\cdots b_{3}b_{1}+\varepsilon_{2k-2}b_{2k-3}b_{2k-5}\cdots b_{3}b_{1}+\varepsilon_{2k-4}b_{2k-5}\cdots b_{3}b_{1}+\cdots+\varepsilon_{2}b_{1}+
    	\varepsilon_{0}+1,\notag \\
    	A(x^{'})&= b_{2k-1}b_{2k-3}\cdots b_{3}b_{1}+(b_{2k-1}-1)b_{2k-3}b_{2k-5}\cdots b_{3}b_{1}+\varepsilon_{2k-4}b_{2k-5}\cdots b_{3}b_{1}+\cdots+\varepsilon_{2}b_{1}+
    	\varepsilon_{0}+1,\notag \\
    	B(x)&=B(x^{'})=b_{2k}b_{2k-2}\cdots b_{2}.\notag
    	\end{align}
    	Then
    	$$\frac{A(x)B(x)}{x}\leq \frac{A(x^{'})B(x^{'})}{x^{'}}\leq \frac{A(y_{k})B(y_{k})}{y_{k}}.$$
    This completes the proof.
    \end{Proof}	
    
     \begin{lemma}
    	For the additive complements $A$ and $B$, if
    	$$x=\varepsilon_{1}a_{1}+\varepsilon_{3}a_{3}+\varepsilon_{5}a_{5}+\cdots+\varepsilon_{2k-1}a_{2k-1}+\varepsilon_{2k+1}a_{2k+1} \in B, 0 \leq \varepsilon_{j}\leq b_{j+1}-1,$$
    	then
    	$$\frac{A(x)B(x)}{x}\leq \frac{A(z_{k})B(z_{k})}{z_{k}}.$$
    \end{lemma}
    
    \begin{Proof}	
    	By noting that $x\in B$, we only need to consider $x\neq z_{k}$, i.e. the following two conditions.
    	
    	\textbf{Case $1$}  $x=\varepsilon_{1}a_{1}+\varepsilon_{3}a_{3}+\varepsilon_{5}a_{5}+\cdots+\varepsilon_{2k-1}a_{2k-1}+\varepsilon_{2k+1}a_{2k+1},2 \leq \varepsilon_{2k+1}\leq b_{2k+2}-1$	.
    	
    	\textbf{Case $2$}  $x=\varepsilon_{1}a_{1}+\varepsilon_{3}a_{3}+\cdots+\varepsilon_{2m-1}a_{2m-1}+\displaystyle \sum_{i=m}^{k-1}(b_{2i+2}-1)a_{2i+1}+a_{2k+1}$, where $\varepsilon_{2m-1}\leq b_{2m}-2,1\leq m \leq k($if $m=k$, then the sum is zero$)$.
    	
    	For Case $1$, we have $A(x)=b_{2k+1}b_{2k-1}\cdots b_{3}b_{1},B(x) \leq (\varepsilon_{2k+1}+1)b_{2k}b_{2k-2}\cdots b_{4} b_{2}.$
    	
    	If $3\leq \varepsilon_{2k+1}\leq b_{2k+2}-1$, we have
    	$$\frac{A(x)B(x)}{x}\leq\frac{(\varepsilon_{2k+1}+1)a_{2k+1}}{\varepsilon_{2k+1}a_{2k+1}}=1+\frac{1}{\varepsilon_{2k+1}}\leq \frac{4}{3} \leq
    	\frac{A(z_{k})B(z_{k})}{z_{k}}.$$
    	If $\varepsilon_{2k+1}=2$, then $x=\varepsilon_{1}+\varepsilon_{3}a_{3}+\varepsilon_{5}a_{5}+\cdots+\varepsilon_{2k-1}a_{2k-1}+2a_{2k+1}$, where $\varepsilon_{2k-1}\leq b_{2k}-1$. Thus $B(x)\leq 2b_{2k}b_{2k-2}\cdots b_{4}b_{2}+(\varepsilon_{2k-1}+1)b_{2k-2}\cdots b_{4}b_{2}$.
    	By Lemma $3.3$ $(x=\varepsilon_{2k-1})$ we have
    	\begin{equation}
    		\begin{aligned}
    			\dfrac{A(x)B(x)}{x}
    			&\leq \dfrac{2a_{2k+1}+(\varepsilon_{2k-1}+1)b_{2k+1}a_{2k-1}}{2a_{2k+1}+\varepsilon_{2k-1}a_{2k-1}}
    			\leq  \dfrac{2a_{2k+1}+(b_{2k}-1+1)b_{2k+1}a_{2k-1}}{2a_{2k+1}+(b_{2k}-1)a_{2k-1}} \\
    			&=\dfrac{3}{2+\frac{1}{b_{2k+1}}-\frac{1}{b_{2k+1}b_{2k}}}
    			\leq \dfrac{2}{1+D_{2k+1}^{*}}
    			=\frac{A(z_{k})B(z_{k})}{z_{k}},\nonumber
    		\end{aligned}
    	\end{equation}
    	where the last inequality is proved by the fact that $0<\frac{1}{b_{2k+1}}-\frac{1}{b_{2k+1}b_{2k}}<\frac{1}{2}$ and $0<D_{2k+1}^{*}-(\frac{1}{b_{2k+1}}-\frac{1}{b_{2k+1}b_{2k}})<\frac{1}{8}$.
    	
    	For Case $2$, if $m\leq k-1$, we have $A(x)=b_{2k+1}b_{2k-1}\cdots b_{3}b_{1}$ and $B(x)\leq b_{2k}b_{2k-2}\cdots b_{4}b_{2}+(b_{2k}-1)b_{2k-2}\cdots b_{4}b_{2}+(b_{2k-2}-1)b_{2k-4}\cdots b_{4}b_{2}+\cdots+(b_{2m+2}-1)b_{2m}\cdots b_{4}b_{2}+(\varepsilon_{2m-1}+1)b_{2m-2}\cdots b_{4}b_{2}$. Hence  by Lemma $3.3$ $(x=\varepsilon_{2m-1})$ we have
    	\begin{equation}
    		\begin{aligned}
    			\dfrac{A(x)B(x)}{x}
    			&\leq \dfrac{\mathscr{M}}
    			{a_{2k+1}+\displaystyle \sum_{i=m}^{k-1}(b_{2i+2}-1)a_{2i+1}+\varepsilon_{2m-1}a_{2m-1}}
    			\leq  \dfrac{\mathscr{M}^{*}}
    			{a_{2k+1}+\displaystyle \sum_{i=m}^{k-1}(b_{2i+2}-1)a_{2i+1}+(b_{2m}-2)a_{2m-1}}\\
    			&=\dfrac{2a_{2k+1}-b_{2k+1}b_{2k-1}b_{2k-3}\cdots b_{2m+1}a_{2m-1}}{a_{2k+1}+\displaystyle \sum_{i=m}^{k-1}(b_{2i+2}-1)a_{2i+1}+(b_{2m}-2)a_{2m-1}}
    			=\dfrac{2-(b_{2k}b_{2k-2}\cdots b_{2m+2}b_{2m})^{-1}}{1+\displaystyle \sum_{i=m}^{k-1}(b_{2i+2}-1)\frac{a_{2i+1}}{a_{2k+1}}+(b_{2m}-2)\frac{a_{2m-1}}{a_{2k+1}} }\\
    			&\leq \dfrac{2}{1+D_{2k+1}^{*}}
    			=\frac{A(z_{k})B(z_{k})}{z_{k}},\nonumber
    		\end{aligned}
    	\end{equation}
    	where $\mathscr{M}=a_{2k+1}+b_{2k+1}(b_{2k}-1)a_{2k-1}+b_{2k+1}b_{2k-1}(b_{2k-2}-1)a_{2k-3}+b_{2k+1}b_{2k-1}b_{2k-2}(b_{2k-3}-1)a_{2k-4}\\
    	+\cdots+b_{2k+1}b_{2k-1}\cdots b_{2m+3}(b_{2m+2}-1)a_{2m+1}+b_{2k+1}b_{2k-1}\cdots b_{2m+1}(\varepsilon_{2m-1}+1)a_{2m-1}$,\\
    	$\mathscr{M}^{*}=a_{2k+1}+b_{2k+1}(b_{2k}-1)a_{2k-1}+b_{2k+1}b_{2k-1}(b_{2k-2}-1)a_{2k-3}+b_{2k+1}b_{2k-1}b_{2k-2}(b_{2k-3}-1)a_{2k-4}+\cdots+\\
    	b_{2k+1}b_{2k-1}\cdots b_{2m+3}(b_{2m+2}-1)a_{2m+1}+b_{2k+1}b_{2k-1}\cdots b_{2m+1}(b_{2m}-2+1)a_{2m-1}.$
    	
    	If $m=k$, then
    	$$x=\varepsilon_{1}a_{1}+\varepsilon_{3}a_{3}+\cdots+\varepsilon_{2k-1}a_{2k-1}+a_{2k+1}, \varepsilon_{2k-1}\leq b_{2k}-2.$$
    	Let
    	$$x^{'}=\varepsilon_{1}a_{1}+\varepsilon_{3}a_{3}+\cdots+(b_{2k}-1)a_{2k-1}+a_{2k+1},$$
    	namely, $\varepsilon_{2k-1}$ in $x$ becomes $b_{2k}-1$ in $x^{'}$ and the rest is the same. We have
    	\begin{align}
    		B(x)&= b_{2k}b_{2k-2}\cdots b_{4}b_{2}+\varepsilon_{2k-1}b_{2k-2}b_{2k-4}\cdots b_{4}b_{2}+\varepsilon_{2k-3}b_{2k-4}\cdots b_{4}b_{2}+\cdots+\varepsilon_{3}b_{2}+
    		\varepsilon_{1}+1,\notag \\
    		B(x^{'})&= b_{2k}b_{2k-2}\cdots b_{4}b_{2}+(b_{2k}-1)b_{2k-2}b_{2k-4}\cdots b_{4}b_{2}+\varepsilon_{2k-3}b_{2k-4}\cdots b_{4}b_{2}+\cdots+\varepsilon_{3}b_{2}+
    		\varepsilon_{1}+1,\notag \\
    		A(x)&=A(x^{'})=b_{2k+1}b_{2k-1}\cdots b_{3}b_{1}.\notag
    	\end{align}
    	Then
    	$$\frac{A(x)B(x)}{x}\leq \frac{A(x^{'})B(x^{'})}{x^{'}}\leq \frac{A(z_{k})B(z_{k})}{z_{k}}.$$   	
      This completes the proof. 
    \end{Proof}	
    
    \textbf{Proof of Lemma $2.1$} \quad For any positive integer $x$, if $x \notin A\cup B$, then
    $$\dfrac{A(x)B(x)}{x}=\dfrac{A(x-1)B(x-1)}{x}<\dfrac{A(x-1)B(x-1)}{x-1}.$$
    So we only select $x \in A\cup B$ to compute $\displaystyle\limsup_{x\rightarrow+\infty}\dfrac{A(x)B(x)}{x}$.
    
    According to Lemma $3.2$, Lemma $3.4$ and Lemma $3.5$, we have 
    $$\displaystyle\limsup_{x\rightarrow+\infty}\dfrac{A(x)B(x)}{x}=
    \displaystyle\limsup_{k\rightarrow+\infty}\frac{2}{1+D_{k}^{*}}.$$
    Since $\displaystyle\liminf_{k\rightarrow+\infty}D_{k}^{*}=\displaystyle\liminf_{k\rightarrow+\infty}D_{k}$, then
    $$\displaystyle\limsup_{x\rightarrow+\infty}\dfrac{A(x)B(x)}{x}=
    \displaystyle\limsup_{k\rightarrow+\infty}\frac{2}{1+D_{k}}.$$ 
    For $x_{k}=a_{2k}-1$, we have $A(x_{k})=b_{1}b_{3}b_{5}\cdots b_{2k-1},B(x_{k})=b_{2}b_{4}b_{6}\cdots b_{2k}$. Thus $A(x_{k})B(x_{k})-x_{k}=1$.
    This completes the proof of Lemma $2.1$.

    \quad \\
    Fang-Yu Ma \\
    School of Mathematics\\
    Shandong University\\    
    Jinan 250100, P.R. China\\   
    E-mail: mfy17864193400@163.com

\end{document}